\documentclass{amsart}

\usepackage{amsmath,amsfonts,amssymb}
\usepackage{mathrsfs}
\usepackage{longtable}

\usepackage{lunadiagrams}

\newtheorem{proposition}{Proposition}[section]

\theoremstyle{definition}
\newtheorem{definition}{Definition}[section]

\numberwithin{equation}{section}

\newcommand{\rank}{\operatorname{rank}}
\newcommand{\C}{\mathbb C}
\newcommand{\R}{\mathbb R}
\newcommand{\Z}{\mathbb Z}

\newcommand{\A}{\mathbf A}
\newcommand{\s}{\mathscr S}

\newcommand{\card}{\mathrm{card}}

\newcommand{\N}{\mathrm N}
\newcommand{\GL}{\mathrm{GL}}
\newcommand{\SL}{\mathrm{SL}}
\newcommand{\SO}{\mathrm{SO}}
\newcommand{\Sp}{\mathrm{Sp}}

\begin{document}

\title{The spherical systems of the wonderful reductive subgroups}
\author{P.\ Bravi and G.\ Pezzini}
\address{Paolo Bravi\\
Dipartimento di Matematica\\
Universit\`a di Roma ``La Sapienza''\\
Piazzale Aldo Moro 5\\
00185, Roma\\
Italy\\
bravi@mat.uniroma1.it}
\address{Guido Pezzini\\
Department Mathematik\\
FAU Erlangen-N\"urnberg\\
Cauerstra\ss e 11\\
91058 Erlangen\\
Germany\\
pezzini@math.fau.de}

\maketitle

\begin{abstract}
We provide the spherical systems of the wonderful reductive subgroups of any reductive group. 
\end{abstract}

\section{Introduction}

In the beginning of the 20th century, real compact connected semi-simple Lie groups $K$, and their involutions $\sigma$, 
were studied and classified. The spaces $K/K^\sigma$, of which a classical example is the sphere $\mathbb S^2=\SO(3)/\SO(2)$,
have the nice property that the multiplicities of irreducible representations of $K$ in the space of continuous functions 
$\mathrm C^0(K/K^\sigma, \R)$ are $\leq 1$.

In \cite{Kr79}, M.~Kr\"amer considered homogeneous spaces $K/L$ (with $L$ closed and connected) having this property, 
suggesting the name {\em spherical} for such subgroups $L$. 
His motivation was that they are useful in the study of $K$-orbits in finite-dimensional $K$-modules. 
He also conjectured that a nice theory of spherical functions could be successfully reproduced on these spaces, 
although this program has been not yet fully carried out. 
Lists of such $K/L$ began to appear in the literature: Kr\"amer himself classified the cases where $K$ is simple, 
I.V.~Mikityuk extended the list in \cite{M87} to include all semi-simple groups, 
and M.~Brion obtained independently in \cite{Br87} the same result, working with complex groups.

Indeed, it is also natural and equivalent to work in the complex setting: 
every compact real Lie group $K$ can be obtained as the group of real points of a unique semi-simple complex group $G$, 
involutions of $K$ extend to involutions of $G$, 
and a reductive algebraic subgroup $H$ of $G$ can be defined to be spherical by requiring that 
the multiplicities of irreducible representations of $G$ in the space of regular functions $\C[G/H]$ be $\leq 1$. 
Such spaces $G/H$ correspond quite closely to weakly symmetric spaces for compact groups, 
introduced by A.~Selberg (see \cite{Se56} and \cite[Main Theorem]{AV99}).

With complex groups, sphericality has also a very useful equivalent characterization: 
one (or equivalently any) Borel subgroup of $G$ has an open orbit in $G/H$. 
This geometrical condition is meaningful for $H$ any algebraic subgroup of $G$, not only for reductive ones. 
Therefore, it was natural to take this as the definition of spherical subgroups of $G$, and try to classify them. 
More in general, a systematic study of {\em spherical varieties}, 
i.e.\ normal $G$-varieties with a dense orbit isomorphic to $G/H$ where $H$ is a spherical subgroup,
began in the '80s (see \cite{LV83} and \cite{BLV86}).

This generalization has been quite fruitful. 
For example, affine smooth spherical varieties are ``local models'' of multiplicity-free hamiltonian manifolds, 
and have been the crucial ingredient to extend the well-known {\em Delzant theorem} from the abelian to the non-abelian case 
(see \cite{D88}, \cite{Lo09b} and \cite{Kn11}).

Moreover, the classification of all spherical subgroups has been completed recently, 
with the contribution of several authors, see \cite{Lu01,Pe04,BP05,B07,Lo09a,BC10,B09,BP11a,BP11b} 
(and also \cite{C09} for a different approach), 
and for this the knowledge of the reductive spherical subgroups has played a key role 
(see for example \cite[Proposition~4.2.3]{Lo09a}).

In this paper, our aim is to clarify the relations between the classical context and this recent more general classification. 
Spherical reductive subgroups of a given $G$ may be infinitely many, 
and we recall how they can all be obtained from a finite list of cases, called {\em spherically closed}, 
by means of a procedure called {\em augmentation}.


Thanks to a fundamental theorem of F.~Knop (see \cite{Kn96}), spherically closed subgroups are {\em wonderful}, i.e.\ $G/H$ admits a particular compactification which provides good invariants of the subgroup $H$. These invariants can be represented by combinatorial objects called {\em spherical systems}, and they encode in an essential way the information on two sets which are naturally associated with the subgroup $H$: the colors and the spherical roots of the homogeneous space $G/H$. The set of colors is the set of prime $B$-divisors in $G/H$. The set of spherical roots is (essentially) the set of equations of the faces of the (simplicial) cone of $G$-invariant valuations on the field of rational functions on $G/H$.

We list these invariants for all reductive wonderful subgroups, and we show some applications. We explain how they can be used to recover easily the $G$-module structure of $\C[G/H]$ for all $H$ reductive and spherical, and how it is possible to characterize the reductive and connected $H$ such that $G/(H,H)$ is also spherical.

We obtain the spherical systems of the reductive wonderful subgroups in a rather indirect way, 
although with a certain effort one could carry out the computation directly. 
Our approach is as follows. 
On one hand, the abstract spherical systems that can be associated with a wonderful reductive subgroup are classified: 
this comes from the combinatorial study of all the abstract spherical systems of \cite{B09} 
and a criterion for reductivity due to Knop \cite{Kn91}. 
On the other hand, reductive spherical subgroups are also classified, as mentioned above; 
therefore, it remains to find out the right matching between the two lists. 
The knowledge of the spherical system associated with $H$ 
prescribes the dimension of $H$ and the rank of its character group: 
this turns out to be enough to uniquely determine the desired correspondence.

Furthermore, this list of the spherical systems of the reductive wonderful subgroups may serve as a definitive reference, 
and be of use for the explicit computation of the spherical subgroup associated with a given spherical system, 
see \cite[Sections~3.1 and~5.2]{BP11a}.

\paragraph{Acknowledgments.}
We thank the anonymous referee for many remarks and suggestions. 
The second-named author was partially supported by the DFG Schwerpunktprogramm 1388 -- Darstellungstheorie.

\section{Generalities}

\paragraph{Notation.}
In the following $G$ will be a semi-simple complex algebraic group of adjoint type, and, unless otherwise stated, 
we will denote by $H$ a wonderful subgroup of $G$ as in Definition~\ref{def:wonderful}.

Let us fix a maximal torus $T$ in $G$ and a Borel subgroup $B$ of $G$ containing $T$. The corresponding set of simple roots of the root system of $(G,T)$ will be denoted by $S$. 

\paragraph{The spherical system associated with a wonderful subgroup.}

\begin{definition}\label{def:wonderful}
A subgroup $H\subseteq G$ is {\em wonderful} if $G/H$ admits a $G$-equivariant compactification $X$ such that: $X$ is complete and non-singular, the open $G$-orbit $G/H\subseteq X$ has complement equal to the union of finitely many prime $G$-divisors, any subset of these prime $G$-divisors has a transversal and non-empty intersection, and these intersections are exactly all $G$-orbit closures of $X$.
\end{definition}

We recall the definition and some properties of certain invariants associated with $G/H$, for $H$ a wonderful subgroup of $G$. For details we refer to \cite{Lu01}.

Let $P_{G/H}$ be the stabilizer of the open $B$-orbit of $G/H$ and denote by $S^p_{G/H}$ the subset of simple roots corresponding to $P_{G/H}$, which is a parabolic subgroup of $G$ containing $B$. 

Let $\Lambda_{G/H}$ be the lattice of $B$-eigenvalues of $B$-eigenvectors in $\mathbb C(G/H)$ and $\mathcal V_{G/H}$ the set of $\mathbb Q$-valued $G$-invariant valuations on $\mathbb C(G/H)$ trivial on the constant functions. Any such valuation $v$ induces the element $\chi\mapsto v(f_\chi)$ of $\mathrm{Hom}(\Lambda_{G/H},\mathbb Q)$, where $f_\chi\in\mathbb C(G/H)$ is a $B$-eigenvector of $B$-eigenvalue $\chi$. Moreover $v$ is uniquely determined by $\chi\mapsto v(f_\chi)$, therefore we may identify $\mathcal V_{G/H}$ with the corresponding subset of $\mathrm{Hom}(\Lambda_{G/H},\mathbb Q)$, which is a convex cone of maximal dimension generated by finitely many elements.

Let $\Sigma_{G/H}\subset\mathbb N\,S$ be the minimal set of primitive elements in $\Lambda_{G/H}$ such that 
\[\mathcal V_{G/H}=\{v\in\mathrm{Hom}(\Lambda_{G/H},\mathbb Q):\langle v,\sigma\rangle\leq0\mbox{ for all }\sigma\in\Sigma_{G/H}\}.\] 
It is called the set of spherical roots of $G/H$, and it is a basis of $\Lambda_{G/H}$. 

We denote by $\Delta_{G/H}$ the set of $B$-stable prime divisors of $G/H$, called colors. A color is said to be moved by a simple root $\alpha$ if it is not stable under the minimal parabolic containing $B$ corresponding to $\alpha$.

For all simple roots $\alpha$ belonging to $\Sigma_{G/H}$
there are exactly two colors moved by $\alpha$, usually denoted by $D_\alpha^+$ and $D_\alpha^-$. In this case set $\A_{G/H}(\alpha)=\{D_\alpha^+, D_\alpha^-\}$. 
Let $\A_{G/H}$ be the union of $\A_{G/H}(\alpha)$ for all $\alpha\in S\cap\Sigma_{G/H}$, 
there is a $\Z$-bilinear pairing $c_{G/H}\colon \Z\A_{G/H}\times\Z\Sigma_{G/H}\to\Z$, called Cartan pairing, 
between colors and spherical roots induced by the valuations of $B$-stable divisors on $\C(G/H)^{(B)}$. 
For all $D\in\A_{G/H}$ and $\sigma\in\Sigma_{G/H}$  we have
\[c_{G/H}(D,\sigma)\leq1.\] 
For all $\alpha\in S\cap\Sigma_{G/H}$ we have also
\[\A_{G/H}(\alpha)=\{D\in\A_{G/H} : c_{G/H}(D,\alpha)=1\},\]
and, for all $\sigma$, the formula
\[c_{G/H}(D_\alpha^+,\sigma)+c_{G/H}(D_\alpha^-,\sigma)=\langle\alpha^\vee,\sigma\rangle\]
holds.

Notice that the objects we have defined are meaningful for any spherical subgroup $H$ of $G$; the only difference is that in general $\Sigma$ is not a basis of $\Lambda$, and it is if and only if $H$ is wonderful. 

The triple $\s_{G/H}=(S^p_{G/H}, \Sigma_{G/H}, \A_{G/H})$ is a Luna spherical system in the sense of \cite[\S2.1]{Lu01}; we will refer to it as the spherical system associated with $H$.

\paragraph{Information from the spherical system.}
The set of colors of $G/H$, considered as an abstract set together with an extension of the Cartan pairing, can be retrieved from $\s_{G/H}$ as the disjoint union $\Delta_{G/H}^a\cup\Delta_{G/H}^{2a}\cup\Delta_{G/H}^b$ where:
\begin{itemize}
\item $\Delta_{G/H}^a=\A_{G/H}$,
\item $\Delta_{G/H}^{2a}=\{D_\alpha : \alpha\in S\cap{\frac1 2}\Sigma_{G/H}\}$,
\item $\Delta_{G/H}^b=\{D_\alpha : \alpha\in S\setminus(S^p_{G/H}\cup\Sigma_{G/H}\cup{\frac1 2}\Sigma_{G/H})\}/\sim$, where $D_\alpha\sim D_\beta$ if $\alpha$ and $\beta$ are orthogonal and $\alpha+\beta\in\Sigma_{G/H}$.
\end{itemize} 
For all $\alpha\in S$ set:
\[\Delta_{G/H}(\alpha)=\left\{\begin{array}{ll}
\varnothing & \mbox{ if $\alpha\in S^p_{G/H}$} \\
\A_{G/H}(\alpha) & \mbox{ if $\alpha\in\Sigma_{G/H}$} \\
\{D_\alpha\} & \mbox{ otherwise}
\end{array}\right.\]
Then $\Delta_{G/H}(\alpha)$ corresponds to the subset of colors that are not stable under the minimal parabolic containing $B$ 
corresponding to $\alpha$.

The full Cartan pairing of $G/H$ is the $\Z$-bilinear map $c_{G/H}\colon\Z\Delta_{G/H}\times\Z\Sigma_{G/H}\to\Z$ such that:
\[c_{G/H}(D,\sigma)=\left\{\begin{array}{ll}
c_{G/H}(D,\sigma) & \mbox{ if $D\in\Delta_{G/H}^a$} \\
{\frac1  2}\langle\alpha^\vee,\sigma\rangle & \mbox{ if $D=D_\alpha\in\Delta_{G/H}^{2a}$} \\
\langle\alpha^\vee,\sigma\rangle & \mbox{ if $D=D_\alpha\in\Delta_{G/H}^b$} 
\end{array}\right.\] 

The dimension of $H$ and the rank of its character group can also be read off the spherical system $\s_{G/H}$. Indeed, the codimension of $H$ in $G$ is equal to 
\[\card\,\Sigma_{G/H}+\dim\,G/P_{G/H},\]
and the rank of its character group $\mathcal X(H)$ is equal to
\begin{equation}\label{eqn:char}
\rank\mathcal X(H) = \card\,\Delta_{G/H}-\card\,\Sigma_{G/H}.
\end{equation}

Finally, the wonderful subgroup $H$ is reductive if and only if there exists $\sigma\in\mathbb N\Sigma_{G/H}$ such that $c_{G/H}(D,\sigma)>0$ for all $D\in\Delta_{G/H}$.

\paragraph{Spherical closure and augmentations.}

Let us consider for a while a spherical subgroup $H$ of a connected reductive complex algebraic $G$, 
without any other assumptions. 
The normalizer $\N_GH$ acts naturally on $G/H$ with $G$-equivariant automorphisms, 
and this induces an action of $\N_GH$ on $\Delta_{G/H}$.

The kernel of this action is a subgroup $\overline H\subseteq \N_GH$, 
which we call the {\em spherical closure} of $H$. If $\overline H=H$, then $H$ is called {\em spherically closed}. 
One can show that $\overline H$ is itself spherically closed (e.g.\ using \cite[Lemma 2.4.2]{BL11}), 
and a theorem of Knop in \cite{Kn96} states that spherically closed subgroups are wonderful
(notice that a spherically closed $H$ always contains the center of $G$).

Thanks to \cite{Lu01}, the subgroup $H\subseteq \overline H$ is uniquely determined by the corresponding {\em augmentation} (of the spherical system $\s_{G/\overline H}$), which is defined as the quintuple $(S^p_{G/\overline H}, \Sigma_{G/\overline H}, \A_{G/\overline H}, \Lambda_{G/H}, c_{G/H})$. More precisely, there is a bijection between all spherical subgroups having spherical closure $\overline H$, and the set of quintuples $(S^p_{G/\overline H}, \Sigma_{G/\overline H}, \A_{G/\overline H}, \Lambda, c)$ satisfying the following conditions:
\begin{enumerate}
\item $\Lambda$ is a lattice of $B$-weights containing $\Sigma_{G/\overline H}$;
\item $c\colon\Z\Delta_{G/\overline H}\times\Lambda\to\Z$ extends the Cartan pairing of $G/\overline H$;
\item for all $\alpha\in\Sigma_{G/\overline H}\cap S$, we have 
$c(D_\alpha^+,\lambda) + c(D_\alpha^-,\lambda) = \langle \alpha^\vee, \lambda\rangle$
for all $\lambda\in\Lambda$;
\item for all $\alpha\in(\frac12\Sigma_{G/\overline H})\cap S$ we have $\alpha\notin\Lambda$, and $\langle\alpha^\vee,\Lambda\rangle\subseteq 2\Z$;
\item for all $\alpha+\beta\in\Sigma_{G/\overline H}$ with $\alpha,\beta\in S$ and $\alpha\perp\beta$, the coroots $\alpha^\vee$ and $\beta^\vee$ coincide on $\Lambda$.
\item for all $\alpha\in S^p_{G/\overline H}$, the coroot $\alpha^\vee$ is zero on $\Lambda$.
\end{enumerate}
Notice that the standard projection $G/H\to G/\overline H$ induces a bijection between the sets of colors of $G/H$ and $G/\overline H$.

It is also easy to deduce from \cite[\S6]{Lu01} that $H_1,H_2\subseteq \overline H$ with $\overline H_1=\overline H_2=\overline H$ satisfy $H_1\subseteq H_2$ if and only if $\Lambda_{G/H_1}\supseteq \Lambda_{G/H_2}$ and $c_{G/H_1}$ restricted to $\Lambda_{G/H_2}$ is equal to $c_{G/H_2}$.

Augmentations can be used to compute easily the {\em weight monoid} of $G/H$ for any spherical subgroup $H$, i.e.\ the set of dominant weights $\Lambda^+_{G/H}$ such that
\[
\C[G/H] = \bigoplus_{\lambda\in\Lambda^+_{G/H}} V(\lambda),
\]
where $V(\lambda)$ is the irreducible $G$-module of highest weight $\lambda$. Weight monoids were already computed in \cite{A10} for the so-called {\em strictly irreducible} affine spherical homogeneous spaces (see {\em loc.cit.} for the definition).

\begin{proposition}
Let $H$ be a spherical subgroup of $G$, and consider the corresponding augmentation $(S^p_{G/\overline H}, \Sigma_{G/\overline H}, \A_{G/\overline H}, \Lambda_{G/H}, c_{G/H})$. Then
\[
\Lambda^+_{G/H} = \left\{ \lambda\in\Lambda_{G/H} \;\middle\vert\; c_{G/H}(D,\lambda)\geq 0 \;,\;\forall D\in\Delta_{G/\overline H}\right\}.
\]
\end{proposition}
\begin{proof}
Any highest weight vector $f\in \C[G/H]$ has no poles, hence his $B$-eigenvalue $\lambda$ is non-negative on any color of $G/H$, with respect to the Cartan pairing of $G/H$.

For the other inclusion, suppose that a $B$-eigenvector $f\in\C(G/H)$ has no pole on the colors of $G/H$. Since the divisor of $f$ is $B$-stable, it follows that $f$ has no pole at all. The variety $G/H$ is smooth, and we conclude that $f$ is a regular function.
\end{proof}

\begin{proposition}
Let $G$ be semisimple and $H\subseteq G$ a connected reductive spherical subgroup. Then its commutator $H'= (H,H)$ is spherical in $G$ if and only if there exists an augmentation $(S^p_{G/\overline H}, \Sigma_{G/\overline H}, \A_{G/\overline H}, \Lambda, c)$ such that $\Lambda\supseteq\Lambda_{G/H}$, $c$ coincides with $c_{G/H}$ on $\Lambda_{G/H}$, and
\begin{equation}\label{eqn:rankcomm}
\rank\Lambda = |\Delta_{G/\overline H}|.
\end{equation}
\end{proposition}
\begin{proof}
Suppose that $H'$ is spherical. Since $H$ and $H'$ are connected the colors of $G/H$ and of $G/H'$ are identified via the standard map $G/H\to G/H'$. The normalizer $\N_GH$ is included in $\N_GH'$, and we obtain the equality $\N_GH=\N_GH'$ because by standard theory of spherical subgroups the quotient $\N_GH'/H'$ is diagonalizable.

We conclude that $\overline H=\overline{H'}$. The subgroup $H'$ corresponds therefore to an augmentation $(S^p_{G/\overline H}, \Sigma_{G/\overline H}, \A_{G/\overline H}, \Lambda_{G/H'}, c_{G/H'})$ satisfying $\Lambda_{G/H'}\supseteq\Lambda_{G/H}$. We easily conclude that $\rank\Lambda_{G/H'} = |\Delta_{G/\overline H}|$ using (\ref{eqn:char}) applied to $\overline H$ and the fact that $\rank(\Lambda_{G/H'}/\Lambda_{G/H})$ is equal to the rank of the character group of $H$.

Suppose now that the augmentation $(S^p_{G/\overline H}, \Sigma_{G/\overline H}, \A_{G/\overline H}, \Lambda, c)$ with the above properties exists. It corresponds to a spherical subgroup $K\subseteq H$ with $\Lambda_{G/K} = \Lambda$, and such that $H/K$ is diagonalizable. As above, we deduce from $\rank\Lambda_{G/K} = |\Delta_{G/\overline H}|$ that $K$ has character group of rank $0$. Hence $K^\circ = H'$, and we obtain that $H'$ is spherical.
\end{proof}

\begin{proposition}\label{prop:normalizer}
Let $H$ be any spherical subgroup of $G$. Then $\N_GH$ is wonderful, and also spherically closed.
\end{proposition}
\begin{proof}
The group $\N_GH$ is wonderful, because thanks to \cite[Corollary 6.3, part a)]{Kn96} the spherical roots of $G/\N_GH$ generate $\Lambda_{G/\N_GH}$. The spherical closure $\overline{\N_GH}$ is also a wonderful subgroup of $G$. We claim that $\N_GH$ and $\overline{\N_GH}$ have the same spherical roots, whence $\Lambda_{G/\N_GH} = \Lambda_{G/\overline{\N_GH}}$, which proves that $\N_GH=\overline{\N_GH}$ thanks to the theory of augmentations recalled above.

We show the equality of the two sets of spherical roots. In \cite[\S7]{Lu01} it is shown that passing from a subgroup $H$ of $G$ to the spherical closure $\overline H$ the set of spherical roots does not change, except for some spherical roots which get doubled. Precisely, a spherical root $\sigma$ is doubled if and only if $\sigma$ is not a simple root, and $2\sigma$ is ``admissible'', in the sense that there exists a wonderful subgroup $L\subset G$ with $\Sigma_{G/L}=\{\sigma\}$ and $S^p_{G/L}=S^p_{G/H}$.

More in general, also when passing from $H$ to $\N_GH$ the spherical roots do not change, except for some of them being doubled.

Suppose now that $\sigma$ is a spherical root of $G/\N_GH$ and gets doubled in passing from $\N_GH$ to $\overline{\N_GH}$. Then either $\frac12\sigma$ or $\sigma$ is a spherical root of $G/H$. Since the ratio of different proportional spherical roots can only be $2$ (regardless of the subgroup of $G$ being considered), we conclude that $\sigma$ is a spherical root of $G/H$.

But the admissibility condition stated above is the same if applied to $H$ or $\N_GH$, therefore $\sigma$ gets doubled already when passing from $H$ to $\overline H$: contradiction.
\end{proof}

\paragraph{Direct product.}
Let $H_1$ and $H_2$ be two subgroups, respectively, of $G_1$ and $G_2$, semi-simple groups of adjoint type. Clearly, $H_1$ and $H_2$ are wonderful (reductive) if and only if their direct product $H_1\times H_2$ is a wonderful (reductive) subgroup of $G_1\times G_2$. In this case, the spherical system 
\[\s_{G_1\times G_2/H_1\times H_2}=((S_1\cup S_2)^p_{G_1\times G_2/H_1\times H_2},\Sigma_{G_1\times G_2/H_1\times H_2}, \A_{G_1\times G_2/H_1\times H_2})\] 
is called direct product of the spherical systems $\s_{G_1/H_1}$ and $\s_{G_2/H_2}$. One has:
\begin{itemize}
\item $(S_1\cup S_2)^p_{G_1\times G_2/H_1\times H_2}=(S_1)^p_{G_1/H_1}\cup(S_2)^p_{G_2/H_2}$,
\item $\Sigma_{G_1\times G_2/H_1\times H_2}=\Sigma_{G_1/H_1}\cup\Sigma_{G_2/H_2}$,
\item $\A_{G_1\times G_2/H_1\times H_2}=\A_{G_1/H_1}\cup\A_{G_2/H_2}$, and for all $D\in\A_{G_i/H_i}$ and $\sigma\in\Sigma_{G_j/H_j}$ the Cartan pairing $c_{G_1\times G_2/H_1\times H_2}(D,\sigma)$ equals $c_{G_i/H_i}(D,\sigma)$ if $i=j$ or vanishes otherwise.
\end{itemize}

\paragraph{Luna diagrams.}
Let us recall how to read the spherical system off its Luna diagram. 

First, the underlying Dynkin diagram is the Dynkin diagram of $G$. 
The subset of vertices of the Dynkin diagram without circles (neither around nor above nor below) corresponds to the set 
$S^p\subseteq S$. 

Each spherical root has its own symbol as in Table~\ref{tab:sr}. 
Equivalence classes of circles (two circles are equivalent if they are joined by a line) 
are in bijective correspondence with colors. Here we precise that a ``zigzag'' line appearing as the fourth case of Table~\ref{tab:sr} (indicating a spherical root $\alpha_1+\ldots+\alpha_r$, where $\alpha_1,\ldots,\alpha_r$ form a subdiagram of type $\mathsf A_r$) is meant in our convention to join the two vertices $\alpha_1$, $\alpha_r$, not the two circles drawn around them. The two circles correspond therefore to two different colors.

A circle corresponds to a color in $\Delta(\alpha)$ if it is around, above or below the vertex corresponding to the simple root 
$\alpha$. 
Let us restrict our attention to the set $\A\subseteq\Delta$. 
The subset of simple roots that are spherical roots corresponds to the subset of vertices which have circles above and below, 
these two circles correspond to the elements of $\A(\alpha)$, $D_\alpha^+$ and $D_\alpha^-$ respectively. 
The Cartan pairing can be retrieved as follows. 
For all $D\in\A$ and all $\sigma\in\Sigma$, $c(D,\sigma)$ equals 1 if and only if $\sigma$ is a simple root such that 
$D\in\A(\sigma)$. 
If $\sigma$ is a spherical root not orthogonal to $\alpha$ then $c(D_\alpha^+,\sigma)$ belongs to $\{1,0,-1\}$ 
and equals $-1$ if and only if there is an arrow starting from the circle corresponding to $D_\alpha^+$ 
and pointing toward $\sigma$. The rest of the Cartan pairing is uniquely determined by the identity 
\[c(D_\alpha^+,\sigma)+c(D_\alpha^-,\sigma)=\langle\alpha^\vee,\sigma\rangle\]
for all $\alpha\in S\cap\Sigma$ and for all $\sigma\in\Sigma$.

\begin{table}\caption{Diagrams of spherical roots}
\begin{center}\label{tab:sr}
\begin{tabular}{cl}
diagram&spherical root\\
\hline
\begin{picture}(600,2100)\put(300,900){\usebox{\aone}}\end{picture}&\begin{picture}(6000,2100)\put(0,600){$\alpha_1$}\end{picture}\\
\begin{picture}(600,2100)\put(300,900){\usebox{\aprime}}\end{picture}&\begin{picture}(6000,1800)\put(0,600){$2\alpha_1$}\end{picture}\\
\begin{picture}(3000,2100)\put(300,0){\line(1,0){2400}}\multiput(0,0)(2400,0){2}{\put(300,0){\line(0,1){600}}\put(300,900){\circle{600}}\put(300,900){\circle*{150}}}\end{picture}&\begin{picture}(6000,1800)\put(0,600){$\alpha_1+\alpha_1'$}\end{picture}\\
\begin{picture}(6000,2100)\put(300,900){\usebox{\mediumam}}\end{picture}&\begin{picture}(6000,1800)\put(0,600){$\alpha_1+\ldots+\alpha_r$}\end{picture}\\
\begin{picture}(3600,2100)\put(0,900){\usebox{\dthree}}\end{picture}&\begin{picture}(6000,1800)\put(0,600){$\alpha_1+2\alpha_2+\alpha_3$}\end{picture}\\
\begin{picture}(7500,2100)\put(300,900){\usebox{\shortbm}}\end{picture}&\begin{picture}(6000,1800)\put(0,600){$\alpha_1+\ldots+\alpha_r$}\end{picture}\\
\begin{picture}(7500,2100)\put(300,900){\usebox{\shortbprimem}}\end{picture}&
\begin{picture}(6000,1800)\put(0,600){$2(\alpha_1+\ldots+\alpha_r)$}\end{picture}\\
\begin{picture}(3900,2100)\put(0,900){\usebox{\bthirdthree}}\end{picture}&\begin{picture}(6000,1800)\put(0,600){$\alpha_1+2\alpha_2+3\alpha_3$}\end{picture}\\
\begin{picture}(9000,2100)\put(0,900){\usebox{\shortcm}}\end{picture}&
\begin{picture}(6000,1800)\put(0,600){$\alpha_1+2(\alpha_2+\ldots+\alpha_{r-1})+\alpha_r$}\end{picture}\\
\begin{picture}(6900,2700)\put(300,1200){\usebox{\shortdm}}\end{picture}&
\begin{picture}(14400,2400)\put(0,900){$2(\alpha_1+\ldots+\alpha_{r-2})+\alpha_{r-1}+\alpha_r$}\end{picture}\\
\begin{picture}(5700,2100)\put(0,900){\usebox{\ffour}}\end{picture}&\begin{picture}(6000,1800)\put(0,600){$\alpha_1+2\alpha_2+3\alpha_3+2\alpha_4$}\end{picture}\\
\begin{picture}(2100,2100)\put(300,900){\usebox{\gsecondtwo}}\end{picture}&\begin{picture}(6000,1800)\put(0,600){$\alpha_1+\alpha_2$}\end{picture}\\
\begin{picture}(2100,2100)\put(300,900){\usebox{\gtwo}}\end{picture}&\begin{picture}(6000,1800)\put(0,600){$2\alpha_1+\alpha_2$}\end{picture}\\
\begin{picture}(2100,2100)\put(300,900){\usebox{\gprimetwo}}\end{picture}&\begin{picture}(6000,1800)\put(0,600){$4\alpha_1+2\alpha_2$}\end{picture}\\
\end{tabular}
\end{center}
\end{table}

\paragraph{Wonderful subgroups.}
We list in the next section the wonderful reductive subgroups $H$ of any adjoint semisimple group $G$, such that $G/H$ cannot be decomposed into a direct product, together with their respective spherical systems, up to automorphisms of $G$.

The correspondence between wonderful reductive subgroups and spherical systems can be deduced indirectly as explained in the introduction. To obtain the list of all wonderful subgroups, we start from the list of connected spherical subgroups of simply connected semisimple groups in \cite{Br87}. For each such $K$, one may check that its normalizer is self-normalizing. In any case $\N_GK$ is wonderful, and even spherically closed, thanks to Proposition~\ref{prop:normalizer}. Then, one must consider other subgroups $H$ of $G$ with $K \subseteq H \subseteq \N_GK$. 

\begin{proposition}
For all cases $H\subseteq G$ listed in Section~\ref{s:list} satisfying $K\subseteq H\subsetneq \N_GK$ (with $K$ as above), the subgroup $H$ of $G$ is wonderful.
\end{proposition}
\begin{proof}
If $G/H$ has rank $1$ then it is wonderful, and if $G/H$ has rank $2$ we refer to \cite{Wa96}. Hence we assume that $G/H$ has rank at least $3$. Checking whether $H$ is wonderful boils down to computing the lattice $\Lambda_{G/H}$, since the spherical roots of $G/H$ are those of $G/\N_GK$ up to multiples.

Notice that for all cases we have $|\N_GK/H|=2$. Then, thanks to \cite[Lemma 2.4]{G11}, to check that each $H$ is wonderful it is enough to check that there exists a spherical root $\sigma\in\Sigma_{\N_GK}$ such that $\frac12\sigma\in\Lambda_{G/H}$. We give the details for the only case where $G$ is of exceptional type, namely case \ref{list:checkwonderfulness}; in the other cases $G$ is classical and the computation is similar.

In case \ref{list:checkwonderfulness} the group $G$ is of type $\mathsf E_7$. Numbering the simple roots $\alpha_1,\ldots,\alpha_7$ as in Bourbaki, the subgroup $K=H$ is a Levi subgroup of the maximal parabolic of $G$ associated with $\alpha_7$.

Checking the list of spherical systems of reductive wonderful subgroups of $\mathsf E_7$ (or directly, using \cite[Lemme 3.1]{Vu90}) one shows that the spherical system of $G/\N_GK$ is that of case \ref{list:checkwonderfulness2}. Its spherical roots are $\tau = 2\omega_1-\omega_6$, $\epsilon = 2\omega_6-2\omega_7-\omega_1$ and $\sigma = 4\omega_7-2\omega_6$, where $\omega_1,\ldots,\omega_7$ are the fundamental dominant weights.

Since $\Lambda_{G/\N_GK}$ has index $2$ in $\Lambda_{G/H}$ and the latter is contained anyway in the set of integral weights, we deduce that $\Lambda_{G/\N_GK}\subset \Lambda_{G/H}\subseteq \Pi = \Z\omega_1\oplus\Z\omega_6\oplus \Z\omega_7$. On the other hand the quotient $\Pi/\Lambda_{G/\N_GK} \cong \Z\oplus\Z\oplus\frac\Z{4\Z}$ has a unique subgroup of order $2$, which implies that the only possibility for $\Lambda_{G/H}$ is $\Lambda_{G/H}= \Z\tau\oplus\Z\epsilon\oplus\Z\frac12\sigma$.
\end{proof}

\section{The list of reductive wonderful subgroups of reductive groups}

\label{s:list}

For notational convenience $G$ will not necessarily be of adjoint type, 
but $H$ will always contain the center $\mathrm Z_G$ of $G$. 

In the exceptional cases $G$ will be implicitly assumed to be simply connected,  
and to describe the connected component of $H\subset G$ we will provide only its semisimple type and fundamental group.

The first case is the trivial one.

\begin{enumerate}
\item\label{list:lasttrivial} $G\subset G$, for all simple (adjoint) groups $G$. The Luna invariants are: $S^p=S$, $\Sigma=\emptyset$, $\A=\emptyset$. 
\end{enumerate}

\paragraph{Symmetric cases.}
We list here those $H$ such that $(G^\sigma)^\circ \subseteq H\subseteq \N_G(G^\sigma)$ for some involution $\sigma\colon G\to G$.
When $G$ is simply connected $G^\sigma$ is connected, hence $H^\circ=G^\sigma$.


Most cases here are very reductive, i.e.\ not contained in any proper parabolic subgroup of $G$. They satisfy $\rank\mathcal X(H)=0$, and are marked in the list with ``VR''. The not very reductive cases are marked with ``He'', since are always of Hermitian type, i.e.\ maximal Levi subgroups, and satisfy $\rank\mathcal X(H)=1$.

For many cases, different values of the parameters produce diagrams that look different. We report them all, and we also do the same if for some values of the parameters the connected component of $H$ is equal to the connected component of some other case. When this happens then the difference in the corresponding diagrams is only in some spherical root that gets doubled.

\begin{enumerate}
\setcounter{enumi}{1}


\item \label{list:GxG} $G\subset G\times G$, for all simple (adjoint) groups $G$, where $G$ is diagonal in $G\times G$ (VR).

Denoting by $\{\alpha_1,\ldots,\alpha_n\}$ and $\{\alpha'_1,\ldots,\alpha'_n\}$ two (equally ordered) copies of the set of simple roots of $G$, we have:
\[
\begin{array}{lcl}
S^p &=& \varnothing,\\
\Sigma &=& \{\alpha_i+\alpha_i': 1\leq i\leq n\}, \\
\A     &=& \varnothing.
\end{array}
\]

\item $\mathrm S(\GL(p)\times\GL(q))\subset \SL(p+q)$, $1\leq p\leq q$ (He).

\noindent
If $p=1$, $q=1$:
\[\begin{picture}(600,1800)\put(300,900){\usebox{\aone}}\end{picture}\]

\noindent
If $p=1$, $q\geq2$:
\[\begin{picture}(6000,1800)\put(300,900){\usebox{\mediumam}}\end{picture}\]

\noindent
If $2\leq p=q$: 
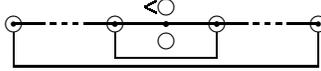
\[\begin{picture}(11400,2400)(-300,-1500)\multiput(0,0)(7200,0){2}{\put(0,0){\usebox{\susp}}\multiput(0,0)(3600,0){2}{\circle{600}}}\multiput(0,-300)(10800,0){2}{\line(0,-1){1200}}\put(0,-1500){\line(1,0){10800}}\multiput(3600,-300)(3600,0){2}{\line(0,-1){900}}\put(3600,-1200){\line(1,0){3600}}\multiput(3600,0)(1800,0){2}{\usebox{\edge}}\put(5400,0){\usebox{\aone}}\put(5400,600){\usebox{\tow}}\end{picture}\]

\noindent
If $1<p<q$, then $\mathrm{card}\,\Sigma=p$, and the diagram is:
\[\diagramaapplusqplusp\]

\item $\N(\mathrm S(\GL(p)\times\GL(p)))\subset \SL(2p)$, $p\geq1$ (VR).

\noindent
If $p=1$:
\[\begin{picture}(600,1800)\put(300,900){\usebox{\aprime}}\end{picture}\]

\noindent
If $p\geq2$:
\[\diagramaaprimepplusoneplusp\]

\item $\SO(p)\cdot\mathrm Z_{\SL(p)}\subset \SL(p)$, $p\geq3$ (VR):
\[\diagramaon\]

\item $\Sp(2p)\cdot\mathrm Z_{\SL(2p)}\subset\SL(2p)$, $p\geq2$ (VR): 
\[\begin{picture}(10800,1800)\put(0,600){\diagramacn}\end{picture}\]

\item $\SO(2p)\subset\SO(2p+1)$, $p\geq2$ (VR):
\[\begin{picture}(7500,1800)\put(300,900){\usebox{\shortbm}}\end{picture}\]

\item $\SO(2)\times\SO(2q-1)\subset\SO(2q+1)$, $q\geq2$ (He).

\noindent
If $q=2$:
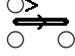
\[\begin{picture}(2400,1800)(-300,-900)
\put(0,0){\usebox{\rightbiedge}}
\put(0,0){\usebox{\aone}}
\put(1800,0){\usebox{\aprime}}
\put(0,600){\usebox{\toe}}
\end{picture}\]

\noindent
If $q\geq 3$:
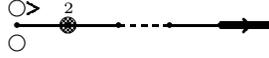
\[\begin{picture}(9300,1800)(-300,-900)
\put(0,0){\usebox{\aone}}
\put(0,0){\usebox{\edge}}
\put(1800,0){\usebox{\shortbprimem}}
\put(0,600){\usebox{\toe}}
\end{picture}\]

\item $\mathrm S(\mathrm O(p)\times\mathrm O(2q+1-p))\subset\SO(2q+1)$, $1\leq p\leq q \geq 2$ (VR).

\noindent
If $p=1$: 
\[\begin{picture}(7500,1800)\put(300,900){\usebox{\shortbprimem}}\end{picture}\]

\noindent
If $p=q=2$:
\[\begin{picture}(2400,1800)(-300,-900)
\put(0,0){\usebox{\rightbiedge}}
\multiput(0,0)(1800,0){2}{\usebox{\aprime}}
\end{picture}\]

\noindent
If $q>p=2$: 
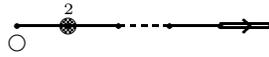
\[\begin{picture}(9300,1800)(-300,-900)
\put(0,0){\usebox{\aprime}}
\put(0,0){\usebox{\edge}}
\put(1800,0){\usebox{\shortbprimem}}
\end{picture}\]

\noindent
If $p=q\geq 3$: 
\[\begin{picture}(9600,1800)(-300,-900)
\multiput(0,0)(5400,0){2}{\usebox{\edge}}
\put(1800,0){\usebox{\susp}}
\put(7200,0){\usebox{\rightbiedge}}
\multiput(0,0)(1800,0){2}{\usebox{\aprime}}
\multiput(5400,0)(1800,0){3}{\usebox{\aprime}}
\end{picture}\]

\noindent
If $2<p<q$, then $\mathrm{card}\,\Sigma=p$, and the diagram is: 
\[\diagrambopplusq\]

\item $\Sp(2p)\times\Sp(2q)\subset\Sp(2p+2q)$, $p\geq1$, $q\geq2$ (VR).

\noindent
If $p=1$: 
\[\begin{picture}(9000,1800)\put(0,900){\usebox{\shortcm}}\end{picture}\]

\noindent
If $p=q$: 
\[\begin{picture}(12900,1800)(0,-900)
\put(0,0){\usebox{\dthree}}
\put(3600,0){\usebox{\susp}}
\put(7200,0){\usebox{\dthree}}
\put(10800,0){\usebox{\leftbiedge}}
\put(12600,0){\usebox{\gcircle}}
\end{picture}\]

\noindent
If $q>p>1$ then $\mathrm{card}\,\Sigma=p$, and the diagram is: 
\[\begin{picture}(19800,1800)(0,-900)
\put(0,0){\usebox{\dthree}}
\put(3600,0){\usebox{\susp}}
\put(7200,0){\usebox{\dthree}}
\put(10800,0){\usebox{\shortcm}}
\end{picture}\]

\item $\N(\Sp(2p)\times\Sp(2p))\subset\Sp(4p)$, $p\geq2$ (VR):
\[\begin{picture}(12900,1800)(0,-900)
\put(0,0){\usebox{\dthree}}
\put(3600,0){\usebox{\susp}}
\put(7200,0){\usebox{\dthree}}
\put(10800,0){\usebox{\leftbiedge}}
\put(12600,0){\usebox{\gcircletwo}}
\end{picture}\]

\item $\GL(p)\subset\Sp(2p)$, $p\geq3$ (He):
\[\begin{picture}(9600,1800)(-300,-900)\multiput(0,0)(5400,0){2}{\usebox{\edge}}\put(1800,0){\usebox{\susp}}\multiput(0,0)(1800,0){2}{\usebox{\aprime}}\multiput(5400,0)(1800,0){2}{\usebox{\aprime}}\put(7200,0){\usebox{\leftbiedge}}\put(9000,0){\usebox{\aone}}\put(9000,600){\usebox{\tow}}\end{picture}\]

\item $\N(\GL(p))\subset\Sp(2p)$, $p\geq3$ (VR):
\[\diagramcoprimen\]

\item $\SO(2)\times\SO(2p)\subset\SO(2p+2)$, $p\geq3$ (He):
\[\begin{picture}(8700,2400)(-300,-1200)
\put(0,0){\usebox{\aone}}
\put(0,0){\usebox{\edge}}
\put(1800,0){\usebox{\shortdm}}
\put(0,600){\usebox{\toe}}
\end{picture}\]

\item $\mathrm S(\mathrm O(p)\times\mathrm O(2q-p))\subset\SO(2q)$, $1\leq p\leq q \geq4$ (VR).

\noindent
If $p=1$: 
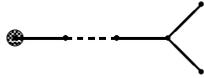
\[\begin{picture}(6900,2400)\put(300,1200){\usebox{\shortdm}}\end{picture}\]

\noindent
If $p=2$: 
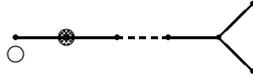
\[\begin{picture}(8700,2400)(-300,-1200)
\put(0,0){\usebox{\aprime}}
\put(0,0){\usebox{\edge}}
\put(1800,0){\usebox{\shortdm}}
\end{picture}\]

\noindent
If $p = q-1$: 
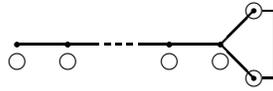
\[\begin{picture}(9450,3000)(-300,-1500)
\multiput(0,0)(5400,0){2}{\usebox{\edge}}
\put(1800,0){\usebox{\susp}}
\put(7200,0){\usebox{\bifurc}}
\multiput(0,0)(1800,0){2}{\usebox{\aprime}}
\multiput(5400,0)(1800,0){2}{\usebox{\aprime}}
\multiput(8400,-1200)(0,2400){2}{\usebox{\wcircle}}
\multiput(8700,-1200)(0,2400){2}{\line(1,0){450}}
\put(9150,-1200){\line(0,1){2400}}
\end{picture}\]

\noindent
If $p=q$: 
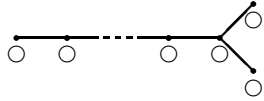
\[\begin{picture}(9000,3600)(-300,-2100)
\multiput(0,0)(5400,0){2}{\usebox{\edge}}
\put(1800,0){\usebox{\susp}}
\put(7200,0){\usebox{\bifurc}}
\multiput(0,0)(1800,0){2}{\usebox{\aprime}}
\multiput(5400,0)(1800,0){2}{\usebox{\aprime}}
\multiput(8400,-1200)(0,2400){2}{\usebox{\aprime}}
\end{picture}\]

\noindent
If $3\leq p\leq q-2$ then $\mathrm{card}\,\Sigma=p$, and the diagram is:
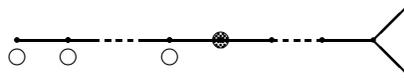
\[\diagramdopplusq\]

\item $\GL(p)\subset\SO(2p)$, $p\geq4$ (He).

\noindent
If $p$ is odd: 
\[\diagramdcnodd\]

\noindent
If $p$ is even: 
\[\begin{picture}(14100,3300)(0,-2100)\multiput(0,0)(7200,0){2}{\usebox{\dthree}}\put(3600,0){\usebox{\susp}}\put(10800,0){\usebox{\edge}}\put(12600,0){\usebox{\bifurc}}\put(12600,0){\usebox{\gcircle}}\put(13800,-1200){\usebox{\aone}}\put(13800,-600){\usebox{\tonw}}\end{picture}\]

\item $\N(\GL(2p))\subset\SO(4p)$, $p\geq2$ (VR):
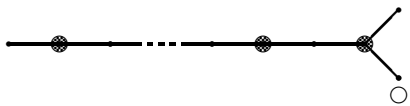
\[\diagramdcprimeneven\]

\item $\mathsf D_5\subset\mathsf E_6$, $H$ maximal Levi subgroup (He):
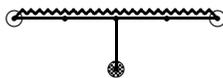
\[\begin{picture}(7800,3000)(-300,-2100)
\put(0,0){\usebox{\afive}}
\put(3600,0){\usebox{\vedge}}
\put(3600,-1800){\usebox{\gcircle}}
\end{picture}\]

\item $\mathsf F_4\subset\mathsf E_6$, $H=H^\circ\cdot\mathrm Z_G$ (VR):
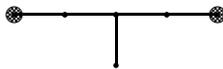
\[\begin{picture}(7800,2700)(-300,-1800)
\multiput(0,0)(1800,0){4}{\usebox{\edge}}
\put(3600,0){\usebox{\vedge}}
\multiput(0,0)(7200,0){2}{\usebox{\gcircle}}
\end{picture}\]

\item $\mathsf A_5\times\mathsf A_1\subset\mathsf E_6$, 
$H$ connected, $\pi_1(H)=\Z/2\Z$ (VR):
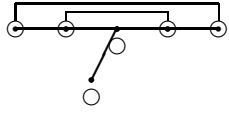
\[\begin{picture}(7800,3600)(-300,-2700)\put(0,0){\usebox{\dynkinafive}}\thicklines\put(3600,0){\line(-1,-2){900}}\thinlines\put(2700,-1800){\circle*{150}}\multiput(0,0)(5400,0){2}{\multiput(0,0)(1800,0){2}{\circle{600}}}\multiput(0,300)(7200,0){2}{\line(0,1){600}}\put(0,900){\line(1,0){7200}}\multiput(1800,300)(3600,0){2}{\line(0,1){300}}\put(1800,600){\line(1,0){3600}}\multiput(3600,0)(-900,-1800){2}{\usebox{\aprime}}\end{picture}\]

\item $\mathsf C_4\subset\mathsf E_6$, 
$H=H^\circ\cdot\mathrm Z_G$, $\pi_1(H^\circ)=\Z/2\Z$ (VR):
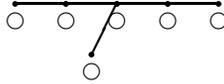
\[\begin{picture}(7800,3600)(-300,-2700)
\put(0,0){\usebox{\dynkinafive}}
\thicklines
\put(3600,0){\line(-1,-2){900}}
\thinlines
\put(2700,-1800){\circle*{150}}
\multiput(0,0)(1800,0){5}{\usebox{\aprime}}
\put(2700,-1800){\usebox{\aprime}}
\end{picture}\]

\item \label{list:checkwonderfulness} $\mathsf E_6\subset\mathsf E_7$, $H$ maximal Levi subgroup (He):
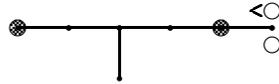
\[\begin{picture}(9600,2700)(-300,-1800)
\multiput(0,0)(1800,0){5}{\usebox{\edge}}
\put(3600,0){\usebox{\vedge}}
\multiput(0,0)(7200,0){2}{\usebox{\gcircle}}
\put(9000,0){\usebox{\aone}}
\put(9000,600){\usebox{\tow}}
\end{picture}\]

\item \label{list:checkwonderfulness2} $\mathsf E_6\subset\mathsf E_7$, $H^\circ$ maximal Levi subgroup, $[H:H^\circ]=2$ (VR):
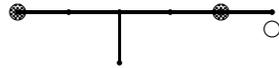
\[\begin{picture}(9600,2700)(-300,-1800)
\multiput(0,0)(1800,0){5}{\usebox{\edge}}
\put(3600,0){\usebox{\vedge}}
\multiput(0,0)(7200,0){2}{\usebox{\gcircle}}
\put(9000,0){\usebox{\aprime}}
\end{picture}\]

\item $\mathsf D_6\times\mathsf A_1\subset\mathsf E_7$, 
$H$ connected, $\pi_1(H)=\Z/2\Z$ (VR):
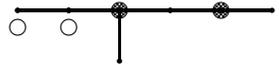
\[\begin{picture}(9600,2700)(-300,-1800)
\multiput(0,0)(1800,0){5}{\usebox{\edge}}
\put(3600,0){\usebox{\vedge}}
\multiput(3600,0)(3600,0){2}{\usebox{\gcircle}}
\multiput(0,0)(1800,0){2}{\usebox{\aprime}}
\end{picture}\]

\item $\mathsf A_7\subset\mathsf E_7$, 
$\pi_1(H^\circ)=\Z/2\Z$, $\mathrm Z_G\subset H^\circ$, $[H:H^\circ]=2$ (VR):
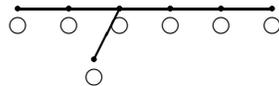
\[\begin{picture}(9600,3600)(-300,-2700)
\put(0,0){\usebox{\dynkinasix}}
\thicklines
\put(3600,0){\line(-1,-2){900}}
\thinlines
\put(2700,-1800){\circle*{150}}
\multiput(0,0)(1800,0){6}{\usebox{\aprime}}
\put(2700,-1800){\usebox{\aprime}}
\end{picture}\]

\item $\mathsf E_7\times\mathsf A_1\subset\mathsf E_8$, 
$H$ connected, $\pi_1(H)=\Z/2\Z$ (VR):
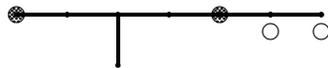
\[\begin{picture}(11400,2700)(-300,-1800)
\multiput(0,0)(1800,0){6}{\usebox{\edge}}
\put(3600,0){\usebox{\vedge}}
\multiput(0,0)(7200,0){2}{\usebox{\gcircle}}
\multiput(9000,0)(1800,0){2}{\usebox{\aprime}}
\end{picture}\]

\item $\mathsf D_8\subset\mathsf E_8$, 
$H$ connected, $\pi_1(H)=\Z/2\Z$ (VR):
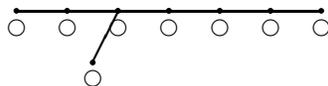
\[\begin{picture}(11400,3600)(-300,-2700)
\put(0,0){\usebox{\dynkinaseven}}
\thicklines
\put(3600,0){\line(-1,-2){900}}
\thinlines
\put(2700,-1800){\circle*{150}}
\multiput(0,0)(1800,0){7}{\usebox{\aprime}}
\put(2700,-1800){\usebox{\aprime}}
\end{picture}\]

\item $\mathsf B_4\subset\mathsf F_4$, $H$ connected and simply connected (VR):
\[\begin{picture}(5700,1800)(0,-900)\put(0,0){\usebox{\ffour}}\end{picture}\]

\item $\mathsf C_3\times\mathsf A_1\subset\mathsf F_4$, 
$H$ connected, $\pi_1(H)=\Z/2\Z$ (VR):
\[\diagramfofour\]

\item $\mathsf A_1\times\mathsf A_1\subset\mathsf G_2$, 
$H$ connected, $\pi_1(H)=\Z/2\Z$ (VR):
\[\begin{picture}(2400,1800)(-300,-900)\put(0,0){\usebox{\lefttriedge}}\multiput(0,0)(1800,0){2}{\usebox{\aprime}}\end{picture}\]

\end{enumerate}

\paragraph{Non-symmetric cases.}
Again, we mark with ``VR'' the very reductive cases, which have $\rank\mathcal X(H)=0$. We do not mark the others, where $\mathcal X(H)$ has always rank $1$. 

\begin{enumerate}
\setcounter{enumi}{30} 

\item $\GL(1)\times\Sp(2p)\subset\SL(2p+1)$, $p\geq2$:
\[\begin{picture}(11400,1800)\put(0,600){\diagramacastn}\end{picture}\] 

\item $\GL(p)\subset\SO(2p+1)$, $p\geq2$:
\[\begin{picture}(11400,1800)(-300,-900)
\put(0,0){\usebox{\atwo}}
\put(1800,0){\usebox{\atwoseq}}
\put(9000,0){\usebox{\btwo}}
\put(10800,0){\usebox{\aone}}
\put(10800,600){\usebox{\tow}}
\end{picture}\]

\item $\N(\GL(p))\subset\SO(2p+1)$, $p\geq2$ (VR):
\[\diagrambcprimen\]

\item $\mathrm{Spin}(7)\subset\SO(9)$ (VR):
\[\begin{picture}(6000,1800)(-300,-900)
\put(0,0){\usebox{\dynkinbfour}}
\multiput(0,0)(5400,0){2}{\usebox{\gcircle}}
\end{picture}\]

\item $\mathsf G_2\subset \SO(7)$ (VR):
\[\begin{picture}(3900,1800)\put(0,900){\usebox{\bthirdthree}}\end{picture}\]

\item $\GL(1)\times\Sp(2p)\subset\Sp(2p+2)$, $p\geq2$:
\[\begin{picture}(9300,1800)(-300,-900)
\put(0,0){\usebox{\shortcm}}
\put(0,0){\usebox{\aone}}
\end{picture}\]

\item $\N(\GL(1)\times\Sp(2p))\subset\Sp(2p+2)$, $p\geq2$ (VR):
\[\begin{picture}(9300,1800)(-300,-900)
\put(0,0){\usebox{\shortcm}}
\put(0,0){\usebox{\aprime}}
\end{picture}\]

\item $\mathsf G_2\cdot \mathrm Z_{\SO(8)}\subset\SO(8)$ (VR):
\[\diagramdsastfour\]

\item $\SO(2)\times\mathrm{Spin}(7)\subset\SO(10)$:
\[\begin{picture}(5400,3000)
\put(300,1500){\usebox{\edge}}
\put(300,1500){\usebox{\aone}}
\put(300,2100){\usebox{\toe}}
\put(1800,0){\diagramdsastfour}
\end{picture}\]

\item $\mathsf A_2\subset\mathsf G_2$, $H$ connected and simply connected (VR):
\[\begin{picture}(2100,1800)(-300,-900)\put(0,0){\usebox{\gtwo}}\end{picture}\]

\item $\mathsf A_2\subset\mathsf G_2$, $H^\circ$ connected and simply connected, $[H:H^\circ]=2$ (VR):
\[\begin{picture}(2100,1800)(-300,-900)\put(0,0){\usebox{\gprimetwo}}\end{picture}\]

\item $\Sp(2)\times\Sp(2p)\times\Sp(2q)\subset\Sp(2p+2)\times\Sp(2q+2)$, $p\geq0,q\geq1$, 
where $\Sp(2)$ is diagonal in $\Sp(2p+2)\times\Sp(2q+2)$ (VR).

Note that $p=q=0$ would also be allowed, it is case (\ref{list:GxG}) for $G$ of type $\mathsf A_1$.

\noindent
If $p=0$:
\[\begin{picture}(12000,1800)(-300,-900)
\put(0,0){\usebox{\vertex}}
\put(0,0){\usebox{\wcircle}}
\put(2700,0){\usebox{\shortcsecondm}}
\multiput(0,-300)(2700,0){2}{\line(0,-1){600}}
\put(0,-900){\line(1,0){2700}}
\end{picture}\]

\noindent
If $p\geq1$:
\[\begin{picture}(17400,1800)(-300,-900)
\multiput(0,0)(9900,0){2}{
\put(0,0){\usebox{\edge}}
\put(1800,0){\usebox{\susp}}
\put(5400,0){\usebox{\leftbiedge}}
\put(1800,0){\usebox{\gcircle}}
}
\multiput(0,0)(9900,0){2}{\usebox{\wcircle}}
\multiput(0,-300)(9900,0){2}{\line(0,-1){600}}
\put(0,-900){\line(1,0){9900}}
\end{picture}\]

\item  $(\Sp(2)\times\Sp(2p)\times\Sp(2q)\times\Sp(2r))\cdot\mathrm Z\subset\Sp(2p+2)\times\Sp(2q+2)\times\Sp(2r+2)$, 
$p,q,r\geq 0$, where $\Sp(2)$ is diagonal in $\Sp(2p+2)\times\Sp(2q+2)\times\Sp(2r+2)$ (VR).

\noindent
If $p=q=r=0$: 
\[\begin{picture}(6000,2700)(-300,-1350)
\multiput(0,0)(2700,0){3}{\usebox{\aone}}
\multiput(0,1350)(2700,0){2}{\line(0,-1){450}}
\put(0,1350){\line(1,0){2700}}
\multiput(0,-1350)(5400,0){2}{\line(0,1){450}}
\put(0,-1350){\line(1,0){5400}}
\multiput(3000,-600)(1050,1200){2}{\line(1,0){1050}}
\put(4050,-600){\line(0,1){1200}}
\end{picture}\]

\noindent
If $p\neq0$ and $q=r=0$: 
\[\begin{picture}(14700,2700)(-300,-1350)
\multiput(0,0)(2700,0){3}{\usebox{\aone}}
\multiput(0,1350)(2700,0){2}{\line(0,-1){450}}
\put(0,1350){\line(1,0){2700}}
\multiput(0,-1350)(5400,0){2}{\line(0,1){450}}
\put(0,-1350){\line(1,0){5400}}
\multiput(3000,-600)(1050,1200){2}{\line(1,0){1050}}
\put(4050,-600){\line(0,1){1200}}
\put(5400,0){\usebox{\shortcm}}
\end{picture}\]

\noindent
If $p,q\neq0$ and $r=0$: 
\[\begin{picture}(20100,2700)(-300,-1350)
\multiput(0,0)(2700,0){2}{\usebox{\aone}}
\put(12600,0){\usebox{\aone}}
\multiput(0,1350)(2700,0){2}{\line(0,-1){450}}
\put(0,1350){\line(1,0){2700}}
\multiput(0,-1350)(12600,0){2}{\line(0,1){450}}
\put(0,-1350){\line(1,0){12600}}
\put(3000,-600){\line(1,0){8250}}
\put(11250,600){\line(1,0){1050}}
\put(11250,-600){\line(0,1){1200}}
\multiput(2700,0)(9900,0){2}{
\put(0,0){\usebox{\edge}}
\put(1800,0){\usebox{\susp}}
\put(5400,0){\usebox{\leftbiedge}}
\put(1800,0){\usebox{\gcircle}}
}
\end{picture}\]

\noindent
If $p,q,r\neq0$:
\[\begin{picture}(27300,2700)(-300,-1350)
\multiput(0,0)(9900,0){3}{\usebox{\aone}}
\multiput(0,1350)(9900,0){2}{\line(0,-1){450}}
\put(0,1350){\line(1,0){9900}}
\multiput(0,-1350)(19800,0){2}{\line(0,1){450}}
\put(0,-1350){\line(1,0){19800}}
\put(10200,-600){\line(1,0){8250}}
\put(18450,600){\line(1,0){1050}}
\put(18450,-600){\line(0,1){1200}}
\multiput(0,0)(9900,0){3}{
\put(0,0){\usebox{\edge}}
\put(1800,0){\usebox{\susp}}
\put(5400,0){\usebox{\leftbiedge}}
\put(1800,0){\usebox{\gcircle}}
}
\end{picture}\]

\item $\GL(1)\times\SL(2)\times\SL(p)\subset\SL(p+2)\times\SL(2)$, $p\geq2$, where $\SL(2)$ is diagonal in $\SL(p+2)\times\SL(2)$.

\noindent
If $p=2$: 
\[\begin{picture}(6900,2700)(-300,-1350)
\multiput(0,0)(3600,0){2}{\usebox{\aone}}
\multiput(0,0)(1800,0){2}{\usebox{\edge}}
\put(1800,0){\usebox{\aone}}
\put(6300,0){\usebox{\aone}}
\multiput(0,1350)(3600,0){2}{\line(0,-1){450}}
\put(0,1350){\line(1,0){3600}}
\multiput(0,-1350)(6300,0){2}{\line(0,1){450}}
\put(0,-1350){\line(1,0){6300}}
\multiput(3900,-600)(1050,1200){2}{\line(1,0){1050}}
\put(4950,-600){\line(0,1){1200}}
\put(0,600){\usebox{\toe}}
\put(3600,600){\usebox{\tow}}
\put(1800,600){\usebox{\tow}}
\end{picture}\]

\noindent
If $p\geq3$:
\[\begin{picture}(10500,2700)(-300,-1350)
\multiput(0,0)(7200,0){2}{\usebox{\aone}}
\multiput(0,0)(5400,0){2}{\usebox{\edge}}
\put(1800,0){\usebox{\shortam}}
\put(9900,0){\usebox{\aone}}
\multiput(0,1350)(7200,0){2}{\line(0,-1){450}}
\put(0,1350){\line(1,0){7200}}
\multiput(0,-1350)(9900,0){2}{\line(0,1){450}}
\put(0,-1350){\line(1,0){9900}}
\multiput(7500,-600)(1050,1200){2}{\line(1,0){1050}}
\put(8550,-600){\line(0,1){1200}}
\put(0,600){\usebox{\toe}}
\put(7200,600){\usebox{\tow}}
\end{picture}\]

\item $\GL(1)\times\SL(2)\times\SL(p)\times\Sp(2q)\subset\SL(p+2)\times\Sp(2q+2)$, $p\geq1$, $q\geq2$, where $\SL(2)$ is diagonal in $\SL(p+2)\times\Sp(2q+2)$.

\noindent
If $p=1$:
\[\begin{picture}(13800,2700)(-300,-1350)
\put(0,0){\usebox{\edge}}
\multiput(0,0)(1800,0){2}{\usebox{\aone}}
\put(4500,0){\usebox{\aone}}
\put(4500,0){\usebox{\shortcm}}
\multiput(0,1350)(4500,0){2}{\line(0,-1){450}}
\put(0,1350){\line(1,0){4500}}
\multiput(1800,-1350)(2700,0){2}{\line(0,1){450}}
\put(1800,-1350){\line(1,0){2700}}
\put(0,600){\usebox{\toe}}
\end{picture}\]

\noindent
If $p=2$:
\[\begin{picture}(15600,2700)(-300,-1350)
\multiput(0,0)(3600,0){2}{\usebox{\aone}}
\multiput(0,0)(1800,0){2}{\usebox{\edge}}
\put(1800,0){\usebox{\aone}}
\put(6300,0){\usebox{\aone}}
\multiput(0,1350)(3600,0){2}{\line(0,-1){450}}
\put(0,1350){\line(1,0){3600}}
\multiput(0,-1350)(6300,0){2}{\line(0,1){450}}
\put(0,-1350){\line(1,0){6300}}
\multiput(3900,-600)(1050,1200){2}{\line(1,0){1050}}
\put(4950,-600){\line(0,1){1200}}
\put(0,600){\usebox{\toe}}
\put(3600,600){\usebox{\tow}}
\put(1800,600){\usebox{\tow}}
\put(6300,0){\usebox{\shortcm}}
\end{picture}\]

\noindent
If $p\geq3$:
\[\begin{picture}(19200,2700)(-300,-1350)
\multiput(0,0)(7200,0){2}{\usebox{\aone}}
\multiput(0,0)(5400,0){2}{\usebox{\edge}}
\put(1800,0){\usebox{\shortam}}
\put(9900,0){\usebox{\aone}}
\multiput(0,1350)(7200,0){2}{\line(0,-1){450}}
\put(0,1350){\line(1,0){7200}}
\multiput(0,-1350)(9900,0){2}{\line(0,1){450}}
\put(0,-1350){\line(1,0){9900}}
\multiput(7500,-600)(1050,1200){2}{\line(1,0){1050}}
\put(8550,-600){\line(0,1){1200}}
\put(0,600){\usebox{\toe}}
\put(7200,600){\usebox{\tow}}
\put(9900,0){\usebox{\shortcm}}
\end{picture}\]

\item\label{list:nsodiag} $\N(\SO(p))\subset\SO(p+1)\times\SO(p)$, diagonally, $4\leq p\leq 6$ (VR).

\noindent
If $p=4$:
\[\begin{picture}(7800,2700)(-300,-1350)
\put(0,0){\usebox{\dynkinbtwo}}
\multiput(0,0)(1800,0){2}{\usebox{\aone}}
\multiput(4500,0)(2700,0){2}{\usebox{\aone}}
\put(0,1350){\line(0,-1){450}}
\multiput(4500,1350)(2700,0){2}{\line(0,-1){450}}
\put(0,1350){\line(1,0){7200}}
\multiput(2100,600)(1050,-1200){2}{\line(1,0){1050}}
\put(3150,600){\line(0,-1){1200}}
\multiput(1800,-1350)(5400,0){2}{\line(0,1){450}}
\put(1800,-1350){\line(1,0){5400}}
\put(0,600){\usebox{\toe}}
\put(1800,600){\usebox{\tow}}
\end{picture}\]

\noindent
If $p=5$:
\[\begin{picture}(8700,2850)(-300,-1500)
\put(0,0){\usebox{\dynkinbtwo}}
\put(4500,0){\usebox{\dynkinathree}}
\multiput(0,0)(1800,0){2}{\usebox{\aone}}
\multiput(4500,0)(1800,0){3}{\usebox{\aone}}
\put(0,1350){\line(0,-1){450}}
\multiput(4500,1350)(3600,0){2}{\line(0,-1){450}}
\put(0,1350){\line(1,0){8100}}
\multiput(2100,600)(1050,-1200){2}{\line(1,0){1050}}
\put(3150,600){\line(0,-1){1200}}
\multiput(0,-1500)(6300,0){2}{\line(0,1){600}}
\put(0,-1500){\line(1,0){6300}}
\multiput(1800,-1200)(6300,0){2}{\line(0,1){300}}
\put(1800,-1200){\line(1,0){4400}}
\put(6400,-1200){\line(1,0){1700}}
\multiput(0,600)(4500,0){2}{\usebox{\toe}}
\multiput(1800,600)(6300,0){2}{\usebox{\tow}}
\end{picture}\]

\noindent
If $p=6$:
\[\begin{picture}(10500,3000)(-2100,-1500)
\put(-1800,0){\usebox{\dynkinbthree}}
\put(4500,0){\usebox{\dynkinathree}}
\multiput(-1800,0)(1800,0){3}{\usebox{\aone}}
\multiput(4500,0)(1800,0){3}{\usebox{\aone}}
\put(0,1500){\line(0,-1){600}}
\multiput(4500,1500)(3600,0){2}{\line(0,-1){600}}
\put(0,1500){\line(1,0){8100}}
\multiput(2100,600)(1050,-1200){2}{\line(1,0){1050}}
\put(3150,600){\line(0,-1){1200}}
\multiput(-1800,1200)(8100,0){2}{\line(0,-1){300}}
\multiput(-1800,1200)(6400,0){2}{\line(1,0){1700}}
\put(100,1200){\line(1,0){4300}}
\multiput(0,-1500)(6300,0){2}{\line(0,1){600}}
\put(0,-1500){\line(1,0){6300}}
\multiput(1800,-1200)(6300,0){2}{\line(0,1){300}}
\put(1800,-1200){\line(1,0){4400}}
\put(6400,-1200){\line(1,0){1700}}
\put(-1800,600){\usebox{\toe}}
\multiput(0,600)(4500,0){2}{\usebox{\toe}}
\multiput(1800,600)(6300,0){2}{\usebox{\tow}}
\end{picture}\]

\item $(\Sp(2)\times\Sp(2)\times\Sp(2p)\times\Sp(2q))\cdot\mathrm Z\subset\Sp(4)\times\Sp(2p+2)\times\Sp(2q+2)$, 
$p\geq0,q\geq1$, where $\Sp(2)\times\Sp(2)$ is embedded in $\Sp(4)$, 
then the first factor $\Sp(2)$ is also embedded into $\Sp(2p+2)$ (complementary to $\Sp(2p)$), 
and the second into $\Sp(2q+2)$ (complementary to $\Sp(2q)$) (VR).


\noindent
If $p=0$:
\[\begin{picture}(16500,2700)(-300,-1350)
\put(0,0){\usebox{\dynkinbtwo}}
\multiput(0,0)(1800,0){2}{\usebox{\aone}}
\multiput(4500,0)(2700,0){2}{\usebox{\aone}}
\put(0,1350){\line(0,-1){450}}
\multiput(4500,1350)(2700,0){2}{\line(0,-1){450}}
\put(0,1350){\line(1,0){7200}}
\multiput(2100,600)(1050,-1200){2}{\line(1,0){1050}}
\put(3150,600){\line(0,-1){1200}}
\multiput(1800,-1350)(5400,0){2}{\line(0,1){450}}
\put(1800,-1350){\line(1,0){5400}}
\put(0,600){\usebox{\toe}}
\put(1800,600){\usebox{\tow}}
\put(7200,0){\usebox{\shortcm}}
\end{picture}\]

\noindent
If $p\geq1$:
\[\begin{picture}(21900,2700)(-300,-1350)
\put(0,0){\usebox{\dynkinbtwo}}
\multiput(0,0)(1800,0){2}{\usebox{\aone}}
\multiput(4500,0)(9900,0){2}{\usebox{\aone}}
\put(0,1350){\line(0,-1){450}}
\multiput(4500,1350)(9900,0){2}{\line(0,-1){450}}
\put(0,1350){\line(1,0){14400}}
\multiput(2100,600)(1050,-1200){2}{\line(1,0){1050}}
\put(3150,600){\line(0,-1){1200}}
\multiput(1800,-1350)(12600,0){2}{\line(0,1){450}}
\put(1800,-1350){\line(1,0){12600}}
\put(0,600){\usebox{\toe}}
\put(1800,600){\usebox{\tow}}
\multiput(4500,0)(9900,0){2}{
\put(0,0){\usebox{\edge}}
\put(1800,0){\usebox{\susp}}
\put(5400,0){\usebox{\leftbiedge}}
\put(1800,0){\usebox{\gcircle}}
}
\end{picture}\]

\item $(\Sp(4)\times\Sp(2p))\cdot \mathrm Z\subset\Sp(4)\times\Sp(2p+4)$, $p\geq1$, where $\Sp(4)$ is diagonal in $\Sp(4)\times\Sp(2p+4)$ (VR).

\noindent
If $p=1$:
\[\begin{picture}(8700,2850)(-300,-1500)
\put(0,0){\usebox{\dynkinbtwo}}
\put(4500,0){\usebox{\dynkincthree}}
\multiput(0,0)(1800,0){2}{\usebox{\aone}}
\multiput(4500,0)(1800,0){3}{\usebox{\aone}}
\put(0,1350){\line(0,-1){450}}
\multiput(4500,1350)(3600,0){2}{\line(0,-1){450}}
\put(0,1350){\line(1,0){8100}}
\multiput(2100,600)(1050,-1200){2}{\line(1,0){1050}}
\put(3150,600){\line(0,-1){1200}}
\multiput(0,-1500)(6300,0){2}{\line(0,1){600}}
\put(0,-1500){\line(1,0){6300}}
\multiput(1800,-1200)(6300,0){2}{\line(0,1){300}}
\put(1800,-1200){\line(1,0){4400}}
\put(6400,-1200){\line(1,0){1700}}
\multiput(0,600)(4500,0){2}{\usebox{\toe}}
\multiput(1800,600)(6300,0){2}{\usebox{\tow}}
\put(6300,600){\usebox{\toe}}
\end{picture}\]

\noindent
If $p\geq2$:
\[\begin{picture}(17400,3000)(-300,-1500)
\put(0,0){\usebox{\dynkinbtwo}}
\put(4500,0){\usebox{\dynkinathree}}
\multiput(0,0)(1800,0){2}{\usebox{\aone}}
\multiput(4500,0)(1800,0){3}{\usebox{\aone}}
\multiput(0,1500)(4500,0){2}{\line(0,-1){600}}
\put(8100,1500){\line(0,-1){50}}
\put(8100,1250){\line(0,-1){350}}
\put(0,1500){\line(1,0){8100}}
\multiput(2100,600)(1050,-1200){2}{\line(1,0){1050}}
\put(3150,600){\line(0,-1){1200}}
\multiput(0,-1500)(6300,0){2}{\line(0,1){600}}
\put(0,-1500){\line(1,0){6300}}
\multiput(1800,-1200)(6300,0){2}{\line(0,1){300}}
\put(1800,-1200){\line(1,0){4400}}
\put(6400,-1200){\line(1,0){1700}}
\multiput(0,600)(4500,0){2}{\usebox{\toe}}
\multiput(1800,600)(6300,0){2}{\usebox{\tow}}
\put(6300,600){\usebox{\tobe}}
\put(8100,0){\usebox{\shortcm}}
\end{picture}\]

\item $\SL(p)\times\GL(1)\subset \SL(p)\times\SL(p+1)$, $p\geq2$, where $\SL(p)$ is diagonal in $\SL(p)\times\SL(p+1)$:
\[\begin{picture}(15900,3900)(-300,-2100)
\put(0,0){\usebox{\edge}}
\put(1800,0){\usebox{\susp}}
\multiput(0,0)(1800,0){2}{\usebox{\aone}}
\put(5400,0){\usebox{\aone}}
\put(5400,600){\usebox{\tow}}
\put(1800,600){\usebox{\tow}}
\put(8100,0){
\put(0,0){\usebox{\edge}}
\put(1800,0){\usebox{\susp}}
\multiput(0,0)(5400,0){2}{\multiput(0,0)(1800,0){2}{\usebox{\aone}}}
\multiput(5400,600)(1800,0){2}{\usebox{\tow}}
\put(1800,600){\usebox{\tow}}
\put(5400,0){\usebox{\edge}}
}
\multiput(0,900)(15300,0){2}{\line(0,1){900}}
\put(0,1800){\line(1,0){15300}}
\multiput(1800,900)(11700,0){2}{\line(0,1){600}}
\put(1800,1500){\line(1,0){11700}}
\multiput(5400,900)(4500,0){2}{\line(0,1){300}}
\put(5400,1200){\line(1,0){4500}}
\multiput(0,-900)(13500,0){2}{\line(0,-1){1200}}
\put(0,-2100){\line(1,0){13500}}
\put(1800,-900){\line(0,-1){900}}
\put(1800,-1800){\line(1,0){9900}}
\multiput(11700,-1800)(0,300){3}{\line(0,1){150}}
\multiput(3600,-1500)(0,300){2}{\line(0,1){150}}
\put(3600,-1500){\line(1,0){6300}}
\put(9900,-1500){\line(0,1){600}}
\multiput(5400,-1200)(2700,0){2}{\line(0,1){300}}
\put(5400,-1200){\line(1,0){2700}}
\end{picture}\]

\item $\N(\SO(p))\subset\SO(p)\times\SO(p+1)$ diagonally, $p\geq7$ (VR).

\noindent
If $p$ is odd:
\[\begin{picture}(17700,5100)(1500,-2400)
\put(1800,0){\usebox{\susp}}
\put(1800,0){\usebox{\aone}}
\put(5400,0){\usebox{\aone}}
\put(11700,0){
\put(0,0){\usebox{\edge}}
\put(1800,0){\usebox{\susp}}
\multiput(0,0)(1800,0){2}{\usebox{\aone}}
\put(5400,0){\usebox{\aone}}
}
\put(5400,0){\usebox{\edge}}
\put(7200,0){\usebox{\rightbiedge}}
\put(17100,0){\usebox{\bifurc}}
\multiput(7200,0)(1800,0){2}{\usebox{\aone}}
\multiput(18300,-1200)(0,2400){2}{\usebox{\aone}}
\put(7200,2700){\line(0,-1){1800}}
\put(7200,2700){\line(1,0){12000}}
\put(19200,2700){\line(0,-1){3300}}
\multiput(19200,1800)(0,-2400){2}{\line(-1,0){600}}
\put(9000,2400){\line(0,-1){1500}}
\put(9000,2400){\line(1,0){9900}}
\put(18900,2400){\line(0,-1){500}}
\put(18900,1700){\line(0,-1){1100}}
\put(18900,600){\line(-1,0){300}}
\multiput(5400,2100)(11700,0){2}{\line(0,-1){1200}}
\put(5400,2100){\line(1,0){1700}}
\put(7300,2100){\line(1,0){1600}}
\put(9100,2100){\line(1,0){8000}}
\multiput(1800,1500)(11700,0){2}{\line(0,-1){600}}
\put(1800,1500){\line(1,0){3500}}
\multiput(5500,1500)(1800,0){2}{\line(1,0){1600}}
\put(9100,1500){\line(1,0){4400}}
\put(9000,-2400){\line(0,1){1500}}
\put(9000,-2400){\line(1,0){9300}}
\put(18300,-2400){\line(0,1){300}}
\multiput(7200,-2100)(9900,0){2}{\line(0,1){1200}}
\put(7200,-2100){\line(1,0){1700}}
\put(9100,-2100){\line(1,0){8000}}
\put(5400,-1800){\line(0,1){900}}
\put(5400,-1800){\line(1,0){1700}}
\put(7300,-1800){\line(1,0){1600}}
\put(9100,-1800){\line(1,0){6200}}
\multiput(15300,-1800)(0,300){3}{\line(0,1){150}}
\multiput(3600,-1500)(0,300){3}{\line(0,1){150}}
\put(3600,-1500){\line(1,0){1700}}
\multiput(5500,-1500)(1800,0){2}{\line(1,0){1600}}
\put(9100,-1500){\line(1,0){4400}}
\put(13500,-1500){\line(0,1){600}}
\multiput(1800,-1200)(9900,0){2}{\line(0,1){300}}
\put(1800,-1200){\line(1,0){1700}}
\multiput(3700,-1200)(1800,0){3}{\line(1,0){1600}}
\put(9100,-1200){\line(1,0){2600}}
\put(1800,600){\usebox{\toe}}
\put(5400,0){\multiput(0,600)(1800,0){2}{\usebox{\toe}}}
\put(9000,600){\usebox{\tow}}
\multiput(13500,600)(3600,0){2}{\usebox{\tow}}
\put(18300,1800){\usebox{\tosw}}
\put(18300,-600){\usebox{\tonw}}
\end{picture}\]

\noindent
If $p$ is even:
\[\begin{picture}(19500,5100)(-300,-2400)
\multiput(0,0)(11700,0){2}{
\put(0,0){\usebox{\edge}}
\put(1800,0){\usebox{\susp}}
\multiput(0,0)(1800,0){2}{\usebox{\aone}}
\put(5400,0){\usebox{\aone}}
}
\put(5400,0){\usebox{\edge}}
\put(7200,0){\usebox{\rightbiedge}}
\put(17100,0){\usebox{\bifurc}}
\multiput(7200,0)(1800,0){2}{\usebox{\aone}}
\multiput(18300,-1200)(0,2400){2}{\usebox{\aone}}
\put(7200,2700){\line(0,-1){1800}}
\put(7200,2700){\line(1,0){12000}}
\put(19200,2700){\line(0,-1){3300}}
\multiput(19200,1800)(0,-2400){2}{\line(-1,0){600}}
\put(9000,2400){\line(0,-1){1500}}
\put(9000,2400){\line(1,0){9900}}
\put(18900,2400){\line(0,-1){500}}
\put(18900,1700){\line(0,-1){1100}}
\put(18900,600){\line(-1,0){300}}
\multiput(5400,2100)(11700,0){2}{\line(0,-1){1200}}
\put(5400,2100){\line(1,0){1700}}
\put(7300,2100){\line(1,0){1600}}
\put(9100,2100){\line(1,0){8000}}
\multiput(1800,1500)(11700,0){2}{\line(0,-1){600}}
\put(1800,1500){\line(1,0){3500}}
\multiput(5500,1500)(1800,0){2}{\line(1,0){1600}}
\put(9100,1500){\line(1,0){4400}}
\multiput(0,1200)(11700,0){2}{\line(0,-1){300}}
\put(0,1200){\line(1,0){1700}}
\put(1900,1200){\line(1,0){3400}}
\multiput(5500,1200)(1800,0){2}{\line(1,0){1600}}
\put(9100,1200){\line(1,0){2600}}
\put(9000,-2400){\line(0,1){1500}}
\put(9000,-2400){\line(1,0){9300}}
\put(18300,-2400){\line(0,1){300}}
\multiput(7200,-2100)(9900,0){2}{\line(0,1){1200}}
\put(7200,-2100){\line(1,0){1700}}
\put(9100,-2100){\line(1,0){8000}}
\put(5400,-1800){\line(0,1){900}}
\put(5400,-1800){\line(1,0){1700}}
\put(7300,-1800){\line(1,0){1600}}
\put(9100,-1800){\line(1,0){6200}}
\multiput(15300,-1800)(0,300){3}{\line(0,1){150}}
\multiput(3600,-1500)(0,300){3}{\line(0,1){150}}
\put(3600,-1500){\line(1,0){1700}}
\multiput(5500,-1500)(1800,0){2}{\line(1,0){1600}}
\put(9100,-1500){\line(1,0){4400}}
\put(13500,-1500){\line(0,1){600}}
\multiput(1800,-1200)(9900,0){2}{\line(0,1){300}}
\put(1800,-1200){\line(1,0){1700}}
\multiput(3700,-1200)(1800,0){3}{\line(1,0){1600}}
\put(9100,-1200){\line(1,0){2600}}
\multiput(0,0)(5400,0){2}{\multiput(0,600)(1800,0){2}{\usebox{\toe}}}
\put(9000,600){\usebox{\tow}}
\multiput(13500,600)(3600,0){2}{\usebox{\tow}}
\put(18300,1800){\usebox{\tosw}}
\put(18300,-600){\usebox{\tonw}}
\end{picture}\]

\end{enumerate}

\end{document}